\documentclass[12pt]{amsart}
\usepackage{amssymb}

\newtheorem{Theorem}{Theorem}[section]
\newtheorem{Lemma}[Theorem]{Lemma}
\newtheorem{Corollary}[Theorem]{Corollary}
\newtheorem{Proposition}[Theorem]{Proposition}

\newtheorem{Definition}[Theorem]{Definition}

\newtheorem{Example}[Theorem]{Example}
\newtheorem{Remark}[Theorem]{Remark}

\def\depth{\mbox{\rm {depth}}}
\def\dim{\mbox{\rm {dim}}}
\def\height{\mbox{\rm {height}}}

\def\Ass{\mbox{\rm {Ass\,}}}
\def\Min{\mbox{\rm {Min\,}}}
\def\max{\mbox{\rm {max\,}}}
\newcommand{\m}{\mathfrak m}

\newcommand{\rar}{\rightarrow}

\def\demo{\noindent{\bf Proof. }}
\def\qed{\hfill$\Box$ \\ \smallskip}

\begin{document}
\title{Depths of Powers of the Edge Ideal of a Tree}

\author{Susan Morey}
\address{Department of Mathematics \\
Texas State University\\
601 University Drive\\ 
San Marcos, TX 78666}
\email{morey@txstate.edu}
\urladdr{http://www.txstate.edu/~sm26/}


\begin{abstract}
Lower bounds are given for the depths of $R/I^t$ for $t \geq 1$ when
$I$ is the edge ideal of a tree or forest.  The bounds are given in
terms of the diameter of the tree, or in case of a forest, the largest
diameter of a connected component and the number of connected
components. These lower bounds provide a lower
bound on the power for which the depths stabilize.
\end{abstract}

\maketitle

\section{Introduction}

Let $G$ be a graph on $n$ vertices and let $R=k[x_1, \ldots , x_n]$ be
a polynomial ring over a field $k$ in $n$ variables. The {\it edge
  ideal} $I=I(G)$ of the graph $G$ is the ideal generated by all
monomials of the form $x_ix_j$ such that $\{x_i, x_j\}$ is an edge of
$G$. Edge ideals of graphs have been
studied by various authors (see for example \cite{Villa},
\cite{SVV}, \cite{AJ}, \cite{B}, and \cite{DFT}). The focus
of this work is to determine the depths of the powers of an edge ideal
of a tree. In particular, a lower bound is given for $\depth(R/I^t)$
when $I$ is the edge ideal of a tree or forest. The bound is given in
terms of the diameter of the tree, or in case of a forest, the largest
diameter of a connected component and the number of connected
components. Note that the lower bounds
on the depths of the ideal of a tree or forest also provide a lower
bound on the power for which the depths stabilize.

Suppose $I$ is the edge ideal of a graph $G$, which is not necessarily a
forest. Then $\depth(R/I)$ has 
been studied by various authors (see for example \cite{GV},
\cite{FHV}, \cite{K}). 
However, relatively little is known about $\depth(R/I^t)$ for
specific values of $t$ other than $t=1$. A notable exception is the
case where all powers of $I$ have a linear resolution, which is
discussed in \cite{HH}. In particular, if the complementary graph of
$G$ is chordal, or $I$ is a square-free Veronese ideal (which includes
the class of complete graphs), then bounds on
$\depth(R/I^t)$ are given in \cite[Section 3]{HH}.

It is known in general that the
depths of the powers of $I$, $\depth(R/I^t)$ stabilize for large
$t$. Indeed this follows from general theorems that apply to any
graded ideal of $R$. In particular, by \cite{Burch}
$ \min \{ \depth(R/I^t)\} \leq n- \ell(I)$ where $\ell(I)$ is the
analytic spread 
of $I$, and the minimum is taken over all powers $t$. In \cite{Brod},
Brodmann showed that for sufficiently large $t$, $\depth(R/I^t)$ is a
constant, and this constant is bounded above by $n-\ell(I)$. It was
shown in \cite[Proposition 3.3]{EH} and an alternate 
proof was given in \cite[Theorem 1.2]{HH} that this is an equality 
for sufficiently large $t$  
when the associated graded ring is Cohen-Macaulay. In general, very
little is known about lower bounds for these depths. 
One partial result is an immediate consequence of \cite[Theorem
  5.9]{SVV}, where it is
shown that if $G$ is bipartite, then $I$ is {\em normally
  torsion-free}, which implies $\Ass(R/I^t)=\Min(R/I)$ for all $t$. 
This implies that the maximal
ideal is not associated to any power of $I$, and thus $\depth(R/I^t)
\geq 1$ for all $t$. In \cite[Proposition 2.1]{HH} it is shown that
$\depth(R/I^t)$ is a nonincreasing function of $t$ when all powers of
$I$ have a linear resolution and conditions are given in that paper
under which all powers of $I$ will have linear quotients. In addition,
in \cite[Corollary 2.6]{HH} a lower bound on $\depth(R/I^t)$ is given
  for a class of ideals satisfying a condition that insures all powers
  have linear quotients.

The main result of this paper is a lower bound for the depth of a
tree, or more generally a forest, given in Theorem ~\ref{forest}:

\smallskip

\noindent {\bf Theorem ~\ref{forest}.} 
Suppose $G$ is a forest with $p$ connected components $G_1,\ldots
,G_p$, and $I=I(G)$. Let $d_i$ be the diameter of
$G_i$, and let $d=\max \{ d_i \} $. Then $\depth(R/I^t)\geq \max \{
\lceil{\frac{d-t+2}{3}}\rceil +p-1,p \}$ for all $t\geq 1$.

\smallskip

\noindent The proof of
the theorem relies on induction, and on using a series of short exact
sequences similar to those used in \cite{HaMorey}. To facilitate the
induction, in Section ~\ref{prelim}, after 
introducing some necessary terminology and notation, a series
of lemmas are proven that determine the depth of a path, and in
Proposition ~\ref{tree} a lower bound based on 
the diameter is given for the depth of any tree. In Section
~\ref{powers} a lower bound is 
first determined on the depths of powers of a path (see Proposition
~\ref{powersOfPaths}), which is then used to prove the main theorem
mentioned above. An additional note in Corollary ~\ref{bonus} provides
an improved lower bound for some trees. Note that these lower bounds
on the depths also provide a lower bound on the power for which the
depths stabilize, as is seen in Remark ~\ref{lowerBound}.


\section{Preliminaries} \label{prelim}

For completeness, some standard notation and terminology from
graph theory and algebra are reviewed here. For additional information,
see \cite{Rafael}. Note that by abuse of notation, $x_i$ will at times be used
to denote both a vertex of a graph $G$ and the corresponding variable
of the polynomial ring. 

Let $G$ be a graph with vertices $\{x_1, \ldots , x_n\}$ and let $x_i$
be a vertex of $G$. The {\it neighbor set} of $x_i$ is $N(x_i)=\{ x_j
\, |\, x_ix_j \, {\mbox{\rm is an edge of}} \, G \}$. A vertex $x_i$ is a
  {\it leaf} if $N(x_i)$ has cardinality one and $x_i$ is {\it
      isolated} if $N(x_i)=\emptyset$. There are two types of operations
    preformed on a graph that produce smaller, related, graphs that
    are referred to as {\it minors} of $G$. The one used here will be
    the {\it deletion}, $G \setminus x$, which is formed by 
removing $x$ from the vertex set of $G$ and deleting any edge in $G$ that
contains $x$. This has the effect of setting $x=0$. 

A special type of
graph that will be used heavily in this paper is a path.

\begin{Definition}
Suppose $n \geq 2$. A {\bf path} $P_n$ of length $n-1$ is a set of $n$
distinct vertices $x_1, \ldots , x_n$ together with $n-1$ edges
$x_ix_{i+1}$ for $1 \leq i \leq n-1$. 
\end{Definition}

The {\it diameter} of a
connected graph is the maximum distance between any two vertices,
where the distance between two vertices is given by the minimum length
of a path connecting the vertices. Thus if the diameter of a graph $G$
is $d$, then there exist vertices $u, v$ of $G$ and a path $P_{d+1}$
of length $d$ connecting $u$ and $v$ such that no path of length less
than $d$ exists between $u$ and $v$. Such a path will be referred to
as a path {\it realizing the diameter} of $G$. 

\medskip

The algebraic notions of analytic spread, associated graded ring, and
associated primes will 
also be needed. For additional background, see \cite{BH} and
\cite{mat}. A prime $P$ is {\it associated} to $I$ if $P=(I:c)$ for   
some $c\in R$. The set of primes associated to $I^t$ will be denoted
$\Ass(R/I^t)$. Note that $\m \in \Ass(R/I^t)$ if and only if
$\depth(R/I^t)=0$, where $\m$ is the homogeneous maximal ideal of
$R$. The set $\Min(R/I)$ consists of all primes that are 
minimal over $I$ with respect to inclusion. In general, $\Min(R/I)
\subseteq \Ass(R/I^n)$ for all $n$. In the case of square-free
monomial ideals, $\Ass(R/I)=\Min(R/I)$ and $I$ is {\it normally
  torsion-free} if and only if $\Ass(R/I^t)=\Min(R/I)$ for all $t$. 

There are three intricately related rings, referred to as {\it blowup
  algebras}, that will be used to determine properties of an ideal
  $I$. The {\it Rees algebra} of an ideal $I$ is the graded subring of
  $R[t]$, where $t$ is an indeterminate, given by 
$$R[It]=R \oplus It \oplus I^2t^2 \oplus \ldots,$$
the {\it associated graded ring} is
$$gr_I(R)=R/I \oplus I/I^2 \oplus I^2/I^3 \oplus \ldots,$$
and the {\it fiber cone} is
$$F(I)=R[It]/\m R[It] = gr_I(R)/\m gr_I(R) = R/\m \oplus I/\m I \oplus
I^2/ \m I^2 \oplus \ldots.$$ 
The {\it analytic spread} $\ell(I)$ is the dimension of the fiber
cone, which is also the minimal number of generators of a minimal
reduction of $I$.

The following basic fact will be useful in several
proofs. For clarity it is stated here.

\begin{Lemma} \label{addVariable}
Let $I$ be an ideal in a polynomial ring $R$, let $x$ be an
indeterminate over $R$, and let $S=R[x]$. Then $\depth S/IS = \depth R/I +1$.
\end{Lemma}

\demo
Note that $S/IS \cong R/I[x] \cong R/I \otimes_k k[x]$ and apply
\cite[Theorem 2.2.21]{Rafael}. 
\qed

Notice that if $x$ is an isolated vertex of a graph $G$, and $G'$ is
the minor formed by deleting $x$, then Lemma ~\ref{addVariable} implies
$\depth(R/I(G)) = \depth(R'/I(G'))+1$ where $R'$ is the polynomial ring
in the variables of $G'$ and $R=R'[x]$. 

There is a well-known result, referred to as the Depth Lemma, that
will be heavily used in the proofs in this paper. The Depth Lemma has
appeared many places in the literature, in multiple forms (see for
example \cite[Lemma 3.1.4]{blowups} or \cite[Proposition 1.2.9]{BH}, also
\cite[Lemma 1.3.9]{Rafael}. Two different versions of the 
lemma will be used in this paper, so both are stated here for ease of
reference. 

\begin{Lemma}$($Depth Lemma$)$ \label{depthLemma}
Let $R$ be a local ring or a Noetherian graded ring with $R_0$
local. If
$$0 \rar A \rar B \rar C \rar 0$$
is a short exact sequence of finitely generated $R-$modules where the
maps are all homogeneous, then $($\cite[Lemma 3.1.4]{blowups}$)$
\begin{itemize}
\item $\depth (A) \geq \depth (B) = \depth (C)$, or
\item $\depth (B) \geq \depth(A) = \depth(C)+1$, or
\item $\depth (C) > \depth(A)=\depth(B)$.
\end{itemize}
Also $($see \cite[Proposition 1.2.9]{BH}$)$
\begin{itemize}
\item $\depth (A) \geq \min \{ \depth (B), \depth (C) +1 \}$,
\item $\depth (B) \geq \min \{ \depth (A), \depth (C)\}$, 
\item $\depth (C) \geq \min \{ \depth (A)-1, \depth (B) \}$.
\end{itemize}
\end{Lemma}

The Depth Lemma will primarily be applied to short exact sequences of
the form given in the lemma below. This type of sequence is well-known
and the proof is elementary and left for the reader. 

\begin{Lemma}\label{exact}
Let $K$ be an ideal of $R$ and let $x$ be an element in $R$. Then the
following sequence is exact: 
$$0 \rar R/(K:x) \stackrel{x}{\rar} R/K \rar R/(K,x) \rar 0.$$
\end{Lemma}

Many of the proofs in this work will use the exact sequence above with
$K$ being a power of $I$ or an ideal, such as $(I^t:y)$, that is related
to a power of $I$. The general technique of using iterated versions of
the sequence above with powers of $I$ and their colons was developed
in \cite{HaMorey}. In particular,
the following result, which appears in the proof of \cite[Theorem
  3.5]{HaMorey}, will prove useful here. Because the result is
contained within the proof and does 
not appear in the statement of the theorem, it is restated here for 
ease of reference. 

\begin{Lemma} \label{RHS}
Let $I$ be a square-free monomial ideal in a polynomial ring $R$ and
let $M$ be a monomial in $R$. If $y$ is a variable such that $y$ does
not divide $M$ and $K$ is the extension in $R$ of the minor of $I$
formed by setting $y=0$, then $((I^t:M),y)=((K^t:M),y)$ for any $t\geq 1$. 
\end{Lemma}

\demo
See the proof of Theorem 3.5 of \cite{HaMorey}.
\qed

While the purpose of this work is to examine lower bounds on the 
depths of $(R/I^t)$, it is useful to note that an obvious upper bound
exists. In general, $\depth(R/I^t) \leq \dim(R/I) = n - \height(I)$. For
trees and graphs in general, this bound can often be mildly
strengthened by knowing that $R/I^t$ is not Cohen-Macaulay, in which
case $\depth(R/I^t) \leq n - \height(I) -1$. For example, it is easy
to show that a path $P_n$ with $n\geq 5$ is never unmixed, and thus in
particular $R/I(P_n)^t$ is not Cohen-Macaulay for $t \geq 1$ and $n
\geq 5$.

When dealing with depths, a lower bound is often needed. A basic
lower bound on the depth of $R/I^t$ exists when $I$ is the edge ideal
of a bipartite graph. 

\begin{Lemma}\label{depth1}
Let $I=I(G)$ for $G$ a bipartite graph. Then $\depth(R/I^t) \geq 1$
for all $t \geq 1$. Moreover, $\depth(R/I^t)=1$ for $t$ sufficiently
large if $G$ is a tree.
\end{Lemma}

\demo
If $G$ is a bipartite graph, then for all $t \geq 1$, all associated
primes of $I^t$ are 
minimal primes of $I$ by \cite[Theorem 5.9]{SVV}. Since the
homogeneous maximal ideal $\m$ is not a minimal prime of $I$ when $I$ is a
square-free monomial ideal, $\m$ is not an associated prime of $R/I^t$
for all $t$. Thus $\depth(R/I^t) \geq 1$ for all $t$.

Now by \cite{SVV} $I$ is normally torsion-free for any bipartite
graph, and so by \cite{Hochster} $R[It]$ is Cohen-Macaulay. Then by
\cite{Huneke} the associated graded ring is Cohen-Macaulay as well. So
by \cite[Proposition 3.3]{EH} or \cite[Theorem 1.2]{HH}, $\depth(R/I^t) =
n-\ell(I)$ for sufficiently large $t$, where $n$ is the number of
vertices of $G$. If $G$ is a tree, then $G$ is of linear type by
\cite[Corollary 3.2]{Villa}, and so is self-reductive. Thus 
$\ell(I)=\nu(I)$ is the minimal number of 
generators of $I$. Since a tree on $n$ vertices has $n-1$ edges,
$\ell(I)=n-1$ and $\depth(R/I^t) =1$ for sufficiently large $t$.
\qed

If $G$ is a graph that is not bipartite, then $G$ contains an odd
cycle, and so by \cite{AJ}, $\m \in \Ass(R/I^t)$ for $t>>0$. Hence
$\depth(R/I^t) =0$ for sufficiently large $t$. In general, this would
not necessarily force equality in Burch's formula, however, if the
graph is not bipartite and has a unique cycle, which is necessarily
odd, equality will hold. Graphs having a 
unique odd cycle satisfy $n=\nu(I)$. By \cite[Corollary 3.2]{Villa}
$I$ is of linear type and thus self-reductive, and $\ell(I)=\nu(I)$. Thus
$n-\ell(I)=0$, forcing equality. 

\medskip

Next the depths of powers of the edge ideal of $P_n$ for
$n$ small are determined. These examples will be used later as the basis for
inductive arguments for more general graphs. 
Unless otherwise specified, when working with $P_n$, the ring $R$ will
be a polynomial ring in $n$ variables over a field. While the notation
$P_n$ actually refers to the path, since the ideal $I=I(P_n)$ is
determined by its monomial generating set, which consists of degree
two monomials corresponding to the edges, by abuse of notation, $P_n$
will also be used to denote this generating set, or the ideal it
generates.

\begin{Example} \label{smallN}
For $n \leq 3$, $\depth(R/I(P_n)^t) =1$ for all $t\geq 1$.
\end{Example}

\demo If $n=1$, then $I=I(P_1)=(0)$, and $\depth (R/I^t)=\depth k[x_1]=1$
  for all $t$.  
If $n=2$, then $I=I(P_2)$ is a complete intersection, and thus $R/I^t$
is Cohen-Macaulay for every power of $t$, hence $\depth(R/I^t)=\dim
  (k[x_1,x_2]/(x_1x_2)^t) =1$ for 
  all $t$ in this case.  

If $n=3$, then $P_3$ has height one and is mixed. Since $\Min (R/I) =
\Min (R/I^t)$ for all $t$ for any monomial ideal, $R/I^t$ is
mixed, and thus not Cohen-Macaulay, for all $t\geq 1$. Since $\dim
R/I^t =2$ for all $t\geq 1$, this implies 
$\depth R/I^t \leq 1$ for all $t$. Now since $P_3$ is bipartite,
combining this with Lemma ~\ref{depth1} yields $\depth
(R/I^t)=1$ for all $t\geq 1$. 
\qed

Let $G$ be a tree or a forest. In order to compute the depth of
$R/I^t$, a bound is first needed for $\depth(R/I)$. As a first case, the
depths of paths will be determined.
Note that since the correspondence between graphs and square-free
monomial ideals of degree two is actually a correspondence between
edges of the graph and generators of the ideal, the ideals in this
paper are primarily considered in terms of their monomial
generating sets. Thus when extending (or contracting) the variables as in Lemma
~\ref{addVariable}, the notation $S/I$ will be used in place of $S/IS$
to simplify notation whenever the generators of $I$ are contained in
the ring $S$.

\begin{Lemma}\label{paths}
If $R=k[x_1, \ldots , x_n]$ and $P=P_n$, then for
$I=I(P_n)$, $\depth (R/I) = \lceil {\frac{n}{3}} \rceil$.
\end{Lemma}

\demo
For $n\leq 3$, this has been shown in Example ~\ref{smallN}. Suppose
$n \geq 4$ and let $I=I(P_n)$. Consider the short exact sequence 
$$0 \rar R/(I : x_{n-1}) \stackrel{x_{n-1}}{\rar} R/I \rar
R/(I,x_{n-1}) \rar 0.$$ 
Now $(I: x_{n-1})=(P_{n-3},x_{n-2},x_n)$, so by induction and Lemma
~\ref{addVariable},
$$\depth(R/(I:x_{n-1})) =  \depth(R'[x_{n-1}]/P_{n-3}) = $$
$$\depth
(R'/P_{n-3}) +1 =
\left\lceil{\frac{n-3}{3}}\right\rceil +1 =
\left\lceil{\frac{n}{3}}\right\rceil$$ 
where $R'=k[x_1, \ldots ,
  x_{n-3}]$. Similarly, 
$(I,x_{n-1})=(P_{n-2}, x_{n-1})$, and by induction and Lemma
~\ref{addVariable},
$$\depth(R/(I,x_{n-1}))=\depth(k[x_1,\ldots,x_{n-2}]/(P_{n-2})) +1 = 
\lceil{\frac{n-2}{3}}\rceil +1.$$ 
Thus by the Depth 
Lemma, since $\depth(R/(I,x_{n-1}))\geq \depth 
(R/(I:x_{n-1}))$, then 
$$\depth(R/I) = \depth
(R/(I:x_{n-1}))=\left\lceil{\frac{n}{3}}\right\rceil.$$  
\qed

The depth formula given above for a path can be extended to a lower
bound for the depth of a tree. Note that since the diameter
is the maximum distance between vertices, a path realizing the
diameter of a tree must connect two leaves of the tree, where a leaf
is a vertex with a unique neighbor. 

\begin{Proposition}\label{tree}
If $G$ is a tree of diameter $d$ and $I=I(G)$, then $\depth(R/I) \geq
\lceil{\frac{d+1}{3}}\rceil$. 
\end{Proposition}

\demo
If $d \leq 2$, then $\lceil{\frac{d+1}{3}}\rceil = 1$, and the
result follows from Lemma ~\ref{depth1}. Thus for $n\leq 3$ the result
holds. Assume $d\geq 3$. Let $u$ and $v$ be vertices of $G$ such that
the distance between $u$ and $v$ is $d$, and let $P_{d+1}$ be a path
connecting $u$ and $v$ that realizes the diameter of $G$. Then $u$ is
a leaf, so let $y$ be the unique neighbor of $u$. Then $(I,y)=(J,y)$
where $J$ is the edge ideal of the minor $G'$ of $G$ formed by
deleting $y$. Notice that the diameter of $G'$ is at least $d-2$ and
that $u$ is idolated in $G'$. Thus if $R'$ is the
polynomial ring formed by deleting $u$,
$\depth(R/(I,y))=\depth(R'[u]/(J,y))=\depth(R'/(J,y))+1 \geq
\lceil{\frac{d-2+1}{3}}\rceil +1\geq \lceil{\frac{d+1}{3}}\rceil$ by
induction and Lemma ~\ref{addVariable}.

Now consider $(I:y)=(K,N(y))$ where $K$ is the ideal of the minor
$G''$ of $G$ formed by deleting the variables in $N(y)$. Let $R''$ be
the polynomial ring formed by deleting the variables in $y \cup
N(y)$. Then the
diameter of $G''$ is at least $d-3$ and $y$ is an isolated vertex, so
$\depth(R/(I:y)) = \depth (R''[y]/K) \geq
\lceil{\frac{d-3+1}{3}}\rceil +1 = \lceil{\frac{d+1}{3}}\rceil$ by
induction and Lemma ~\ref{addVariable}. 

The result now follows from applying the Depth Lemma to the sequence
$$0 \rar R/(I:y) \rar R/I \rar R/(I,y)\rar 0.$$
\qed

The goal of this paper is to examine the depths of powers of ideals. A
final preliminary lemma is needed to facilitate calculating the depths
of powers through induction on the power.

\begin{Lemma}\label{leaf}
Suppose $G$ is a graph, $I=I(G)$, $x$ is a leaf of $G$, and $y$ is the
unique neighbor of $x$. Then $(I^t:xy)=I^{t-1}$ for any $t \geq 2$. 
\end{Lemma}

\demo
Since $\{x,y\}$ is an edge of $G$, $xy$ is a generator of $I$ and one
inclusion is clear. Now 
let $a$ be a monomial generator of $(I^t:xy)$. Then
$axy = e_1 \cdots e_t h$ for some degree two monomials $e_i$
corresponding to edges of $G$ and some 
monomial $h$. If $a \not\in I^{t-1}$, then $x$ divides $e_j$ and $y$
divides $e_k$ for
some $j\not=k$. We may assume $j=t$. But since $x$ is a leaf of $G$,
$e_t=xy$ and thus $a=e_1 \cdots e_{t-1} h \in I^{t-1}$.
\qed

\begin{Corollary}\label{powersReduce}
For $n\geq 2$ and $t \geq 2$,
$$(P_n^t : x_{n-1}x_n )= P_n^{t-1}.$$
\end{Corollary}

\demo
Notice that $x_n$ is a leaf of $P_n$ and apply Lemma ~\ref{leaf}.
\qed


\section{Powers of Trees and Forests}\label{powers}

The goal of this paper is to use graph invariants to provide lower
bounds on the depths of the powers of the edge ideal of a tree.
When the graph is a tree or forest, a lower bound on the depth of any
power will be given in Theorem ~\ref{forest}. Since the proof makes
repeated use of applying the 
Depth Lemma to a pair of
sequences, we first prove a lemma to simplify the main proof.

\begin{Lemma}\label{sequencePair}
Suppose $I=I(G)$ for a graph $G$, $z_1$ and $z_2$ are vertices of $G$,
and for some $s\geq 0$, $\depth(R/(I^t:z_1z_2))\geq s$, $\depth(R/(I^t,z_1))\geq s$,
and $\depth(R/((I^t:z_1),z_2)) \geq s$,
then $\depth(R/I^t) \geq s$. 
\end{Lemma} 

\demo
Applying the Depth Lemma to the short exact sequence
$$0 \rar R/(I^t:z_1z_2) {\stackrel{\cdot z_2}{\rar}} R/(I^t:z_1) \rar
R/((I^t:z_1),z_2) \rar 0$$
yields $\depth
(R/(I^t:z_1)) \geq s$. Now apply the Depth Lemma a second time
to the sequence
$$0 \rar R/(I^t:z_1) {\stackrel{\cdot z_1}{\rar}} R/I^t \rar
  R/(I^t,z_1) \rar 0$$ 
to see that $\depth(R/I^t) \geq s$.
\qed

As a first step toward determining the depth of powers of edge ideals
of trees and forests, we can now determine a lower bound on the depth
of the powers of a path ideal. 

\begin{Proposition} \label{powersOfPaths}
For $P_n$ a path ideal with $n\geq 2$, $\depth (R/P_n^t) \geq \max \{ \lceil
{\frac{n-t+1}{3}} \rceil ,1\}$. 
\end{Proposition}

\demo
Notice that since $P_n$ is a bipartite graph, $\depth(R/P_n^t) \geq 1$
for all $t$ by \cite{SVV}, as seen in Lemma ~\ref{depth1}. Thus the
focus of the proof is to show that 
$\depth(R/P_n^t) \geq \lceil{\frac{n-t+1}{3}}\rceil$. The proof is by
induction on $n$ and $t$. Notice that by Example 
~\ref{smallN} the result holds for $n \leq 3$ for all $t$, and by Lemma
~\ref{paths} the result holds for $t=1$ for all $n$. Assume $n\geq 4$
and $t \geq 2$.
Notice that $(P_n^t, x_{n-1}) = (P_{n-2}^t,x_{n-1})$ since $x_{n-1}$ is
the unique neighbor of $x_n$. By induction on $n$, $\depth
R''/P_{n-2}^t \geq \lceil {\frac{n-2-t+1}{3}} \rceil $ where $R''$ is the
polynomial ring in $n-2$ variables. Thus 
$$\depth (R/(P_n^t,x_{n-1}))= \depth
(R''[x_{n-1},x_n]/(P_{n-2}^t,x_{n-1})) =$$
$$\depth(R''/P_{n-2}^t)+1 
\geq \left\lceil {\frac{n-2-t+1}{3}} \right\rceil +1 = \left\lceil {\frac{n-t+2}{3}}
 \right\rceil .$$ 

By Corollary ~\ref{powersReduce} and induction on $t$, 
$$\depth (R/(P_n^t : x_{n-1}x_n)) = \depth (R/P_n^{t-1}) \geq $$
$$\left\lceil {\frac{n-(t-1)+1}{3}}\right\rceil = \left\lceil
       {\frac{n-t+2}{3}} \right\rceil .$$ 

To find the depth of $((P_n^t:x_{n-1}),x_n)$,
note that since 
$x_n$ does not divide $x_{n-1}$,
$((P_n^t:x_{n-1}),x_n)=((P_{n-1}^t:x_{n-1}),x_n)$ by Lemma
~\ref{RHS}. Let $R'=k[x_1,\ldots , x_{n-1}]$ and notice that $\depth
(R/((P_n^t:x_{n-1}),x_n))=\depth
(R'/(P_{n-1}^t:x_{n-1}))$. Consider the
short exact sequence
$$0 \rar R'/(P_{n-1}^t : x_{n-1}x_{n-2}) \rar R'/(P_{n-1}^t:x_{n-1}) \rar R'/((P_{n-1}^t :x_{n-1}),
x_{n-2}) \rar 0.$$ 
By Corollary ~\ref{powersReduce}, $(P_{n-1}^t :
x_{n-1}x_{n-2})=P_{n-1}^{t-1}$, so by induction, 
$$\depth  R'/(P_{n-1}^t
: x_{n-1}x_{n-2}) \geq \left\lceil {\frac{n-1 - (t-1) +1}{3}}\right\rceil =
\left\lceil {\frac{n-t+1}{3}}\right\rceil.$$
 Also, by Lemma ~\ref{RHS}, $((P_{n-1}^t :x_{n-1}),
x_{n-2})=((P_{n-3}^t : x_{n-1}),x_{n-2})=(P_{n-3}^t, x_{n-2})$.  Now
by induction on $n$,
$$\depth(R'/(P_{n-3}^t,x_{n-2})) = \depth (k[x_1, \ldots
  x_{n-3},x_{n-1}]/P_{n-3}^t) =$$
$$\depth (k[x_1,\ldots x_{n-3}]/P_{n-3}^t)
+1 \geq \left\lceil {\frac{n-3 -t+1}{3}}\right\rceil +1 = \left\lceil
{\frac{n-t+1}{3}} \right\rceil.$$
By applying the Depth
Lemma to the sequence above, 
$\depth (R'/(P_{n-1}^t :x_{n-1})) \geq \lceil {\frac{n-t+1}{3}}
\rceil $, so $\depth (R/((P_n^t:x_{n-1}),x_n)) \geq
\lceil{\frac{n-t+1}{3}} \rceil$. 
The result follows from Lemma ~\ref{sequencePair}. 
\qed

The final lemma is an elementary result about trees that will be needed
in the proof of the theorem. 

\begin{Lemma}\label{lotsOfLeaves}
If $G$ is a tree and $P_{d+1}=\{x_1x_2,x_2x_3,\ldots ,x_dx_{d+1}\}$ is
a path realizing 
the diameter of $G$, then at most one element of $N(x_d)$ is not a leaf.
\end{Lemma}

\demo
Note that $P$ is a path of maximal length in $G$
since for a tree, there is a unique path connecting any two
vertices. Let $x \in N(x_d)$. If $x \not= x_{d-1}$, and $x$ is not a 
leaf, then there exists a vertex $z\in N(x)$, $z \not= x_d$. Then the
path $P=\{ x_1x_2,\ldots , x_{d-1}x_d,x_dx,xz\}$ has length $d+1$, which is
a contradiction to $d$ being the diameter of $G$. Thus at
most one neighbor of $x_d$ is not a leaf.
\qed

We are now ready to prove the main theorem. To
simplify the wording, the phrase {\it connected component} of $G$ will
refer only to components containing 
at least two vertices. Isolated vertices will not be considered as
connected components of $G$.

\begin{Theorem}\label{forest}
Suppose $G$ is a forest with $p$ connected components $G_1,\ldots
,G_p$, and $I=I(G)$. Let $d_i$ be the diameter of
$G_i$, and let $d=\max \{ d_i \} $. Then $\depth(R/I^t)\geq \max \{
\lceil{\frac{d-t+2}{3}}\rceil +p-1,p \}$ for all $t\geq 1$.
\end{Theorem}

\demo
The proof is by induction on $t$ and on $n$ where $n$ is the number of
non-isolated vertices of $G$. Without loss of generality, assume
$d=d_1$. 
For $t=1$ and $p=1$ the result follows from Proposition ~\ref{tree}. For
$t=1$ and any $p 
\geq 2$ the result follows from \cite[Lemma 6.2.7]{Rafael} and Proposition
~\ref{tree}. Thus the result 
holds for $t=1$ for any value of $n$. 
Assume $t\geq 2$. 
If $n=2$, then $G=P_2$, $d=1$, and the result holds for all
$t$ by Proposition
~\ref{powersOfPaths}. Assume $n\geq 3$.
Fix a path $P_{d+1}$ in $G$ realizing the diameter,
let $x_1$ be an endpoint of this path (and thus a leaf of $G$), let
$y$ be its unique 
neighbor, and let $N(y)=\{x_1,\ldots ,x_r\}$ be the neighbors of
$y$. Note that $r\geq 1$ and $r$ is finite. Note also that by Lemma
~\ref{lotsOfLeaves}, at most one $x_i$ is not a leaf. Without loss of
generality, assume $x_i$ is a leaf for $1\leq i <r$. Let $I_j$ be the ideal of
the minor of $G$ formed by deleting $x_1,\ldots ,x_j$ for any $1 \leq
j \leq r$. Let $R_j=k[x_{j+1}, \ldots ,x_{n-1},y]$ be the subring of
$R$ excluding $x_1, \ldots, x_j$, and let
$R_j'=k[x_{j+1},\ldots ,x_{n-1}]$. Notice that for each $j$,
$I_j\subset R_j$ is
the edge ideal of a graph involving fewer than $n$ vertices.  

To use Lemma ~\ref{sequencePair} to find the depth of $R/I^t$, the
depths of three ideals must be checked. For ease of notation, let
$s=\max \{ \lceil{\frac{d-t+2}{3}}\rceil +p-1,p \}$.    
Since $x_1$ is a leaf, by Lemma ~\ref{leaf},
$(I^t:x_1y)=I^{t-1}$, and so by induction on $t$,
$\depth(R/(I^t:x_1y))\geq \max \{
\lceil{\frac{d-(t-1)+2}{3}}\rceil +p-1,p\}\geq s$. 

\medskip

To find the depth of the second ideal, note that $(I^t,y)=(J^t,y)$
where $J$ is the edge ideal of the minor $G'$ of $G$ 
formed by deleting $y$. Then $G'$ is again a forest with fewer than
$n$ vertices, at least $p-1$ connected components, and the generators
of $J$ live in $R_1'$. Thus by induction, $\depth(R_1'/J^t)\geq p-1$, so 
$$\depth(R/(I^t,y)) = \depth(R_1[x_1]/(J^t,y)) \geq p-1 + 1 = p.$$
 Suppose
$d\leq 3$. Then $\lceil{\frac{d-t+2}{3}}\rceil \leq 1$ for all $t \geq 2$,
  and thus $s=p$, and $\depth(R/(I^t,y)) \geq s$.   

For $d>3$, note that the number of connected
  components of $J$ is at least $p$ since $G_2, \ldots, G_p$ and $d-2
  \geq 1$ edges of $P_{d+1}$ survive in $G'$. This also implies that
  the maximal diameter $d'$ of a connected component of $G'$ is at 
  least $d-2$. Thus by induction on $n$, 
$$\depth(R_1[x_1]/(J^t,y)) \geq
  \max \left\{\left\lceil{\frac{d-2-t+2}{3}}\right\rceil
  +p-1,p\right\} +1 \geq s.$$  
Thus for all $d$, $\depth (R/(I^t,y)) \geq s$.

\medskip

Now consider $((I^t:y),x_1)$. By Lemma ~\ref{RHS},
$((I^t:y),x_1)=((I_1^t:y),x_1)$ where $I_1$ is as defined above. Thus
$\depth(R/((I^t:y),x_1))=\depth(R_1/(I_1^t:y))$.
Note that if $r=1$, then $d=1$ and $s=p$ as above. Also, for $r=1$,
$(I_1^t:y)=I_1^t$, 
and $I_1\subset R_1'$ is the edge ideal of a forest with fewer than
$n$ variables 
and $p-1$ connected components. Thus
$\depth(R/((I^t:y),x_1))=\depth(R_1/(I_1^t))=\depth(R_1'[y]/I_1^t)
\geq p-1+1=p$. 
Thus for $r=1$, $\depth(R/((I^t:y),x_1)) \geq s$.

Suppose $r \geq 2$. Then $d \geq 2$ since $x_1,y,x_r$ are all vertices
included in $P_{d+1}$. To find $\depth(R_1/(I_1^t:y))$, use reverse
induction on $r$. First consider
$(I_r^t:y)=I_r^t$. Notice that the generators of $I_r^t$ all lie in
$R_r'$, so
$\depth(R_r/I_r^t)=\depth(R_r'[y]/I_r^t)=\depth(R_r'/I_r^t) +1$. 
If $d \leq 3$, then $I_r\subset R_r'$ corresponds to a  graph with
$p-1$ connected components. So by induction on $n$,
$\depth(R_r'/I_r^t) \geq p-1$, so $\depth(R_r/I_r^t) \geq p=s$. 
If $d \geq 4$, $I_r$ corresponds to a graph with diameter at
least $d-3$ with $p$ connected components, so by induction on $n$,
$\depth(R_r'/I_r^t) \geq \max 
\{ \lceil{\frac{d-3-t+2}{3}}\rceil +p-1, p \}$, and thus
$\depth(R_r/(I_r^t:y))=\depth(R_r/I_r^t)=\depth(R_r'/I_r^t) +1 \geq s$. 

Now assume $\depth(R_j/(I_j^t:y))\geq s$ for some $2 \leq j \leq
r$. Consider the short exact sequence
$$0 \rar R_{j-1}/(I_{j-1}^t:yx_j) \rar R_{j-1}/(I^t_{j-1}:y) \rar
R_{j-1}/((I_{j-1}^t:y),x_j)\rar 0.$$
Then by assumption, since $((I_{j-1}^t:y),x_j)=((I_j^t:y),x_j)$ by
Lemma ~\ref{RHS}, 
$\depth(R_{j-1}/((I_{j-1}^t:y),x_j))=\depth(R_j/(I_j^t:y)) \geq s$. 
By Lemma ~\ref{leaf}, $(I_{j-1}^t:yx_j)=I_{j-1}^{t-1}$ since $x_j$ is a
leaf for $j<r$ and $y$ is a leaf when $j=r$. The diameter of the graph
associated to 
$I_{j-1}$ is at least $d-1$, and so 
$$\depth(R_{j-1}/(I_{j-1}^t:yx_j))=\depth(R_{j-1}/I_{j-1}^{t-1}) \geq$$
$$\max \left\{ \left\lceil{\frac{d-1-(t-1)+2}{3}}\right\rceil +p-1, p
 \right\} =s$$
 by
induction. Thus by
the Depth Lemma, $\depth (R_{j-1}/(I^t_{j-1}:y)) \geq s$. Hence by
reverse induction, $\depth(R_1/(I_1^t:y))\geq s$. 

Since
$\depth(R/((I^t:y),x_1))=\depth(R_1/(I_1^t:y))$, then as above
$\depth(R/((I^t:y),x_1)) \geq s$. Thus by applying Lemma
~\ref{sequencePair} with $z_1=y$ and $z_2=x_1$, $\depth(R/I^t) \geq s$
as desired.
\qed

\begin{Corollary}\label{treePowers}
If $G$ is a tree of diameter $d$ and $I=I(G)$, then $\depth(R/I^t) \geq \max \{
\lceil{\frac{d-t+2}{3}}\rceil, 1\}$ for all $t\geq 1$.
\end{Corollary}

\demo
Apply Theorem ~\ref{forest} with $p=1$.
\qed

Note that a path $P_n$ has diameter $d=n-1$, so Corollary
~\ref{treePowers} agrees with Lemma ~\ref{powersOfPaths} for this
special case. The proof above depends heavily on the existance of a
vertex $y$ at most one of whose neighbors is not a leaf. A careful
examination of Lemma ~\ref{lotsOfLeaves} guarantees that any tree with
diameter $d \geq 3$ will
contain at least two such vertices that are not themselves leaves,
namely the neighbors of the two 
leaves of a path realizing the diameter. Call a vertex $v$ of $G$ a
{\it near leaf} of $G$ if $v$ is not a leaf and $N(v)$ contains at
most one vertex that is not a leaf. Let $q$ denote the 
number of near leaves of $G$. Then the bound given
in Theorem ~\ref{forest} can be strengthened using essentially the same
proof. However, a strengthening of Proposition ~\ref{tree} is
needed.

\begin{Lemma}\label{treeBonus}
If $G$ is a tree of diameter $d\geq 1$, $I=I(G)$, and $G$ has $q$ near
leaves, then $\depth(R/I) \geq \lceil{\frac{d+q-1}{3}}\rceil$. 
\end{Lemma}

\demo
For small values of $n$, $q \leq 2$ and the result
holds by Proposition ~\ref{tree}, so
assume $q \geq 3$. Note that for a connected graph, if two 
near leaves are adjacent, $d=3$ and $q=2$ since all other vertices
must be leaves. Thus for $q \geq 3$, no neighbor of a near leaf is a
near leaf. Let $P_{d+1}$ be a path realizing the diameter of
$G$ with vertices $x_1,x_2, \ldots x_{d+1}$. Note that $x_2$ and $x_d$
are both near leaves, and $d \geq 4$ since for $d=3$, $x_2$ and $x_d$
are adjacent, and for $d \leq 2$ a tree has at most one near leaf.

Consider $(I,x_2)=(J,x_2)$ where $J$ is the ideal of the minor $G'$ of $G$
formed by deleting $x_2$. Let $R'$ be the polynomial ring formed by
deleting $x_1$ and $x_2$. The diameter of $G'$ is at least $d-2$ and 
$G'$ has at least $q-1$ near leaves. Thus by induction
$\depth(R'/J) \geq \lceil{\frac{d-2+q-1-1}{3}}\rceil$. Thus
by Lemma ~\ref{addVariable}, $\depth(R/(I,x_2))=\depth(R'[x_1]/J)\geq
\lceil{\frac{d+q-3-1}{3}}\rceil +1\geq \lceil{\frac{d+q-1}{3}}\rceil$. 

Now consider $(I:x_2)=(K,N(x_2))$ where $K$ is the ideal of the minor
$G''$ of $G$ formed by deleting the vertices in $N(x_2)$. The diameter
of $G''$ is at least $d-3$. Let $a$ denote the number of near leaves
adjacent to $x_3$ but not on $P_{d+1}$. Note that any path from 
$x_3$ to a leaf where the path does not contain $x_i$ for $i\not=3$
must have length at most two, else there exists a path of length
greater than $d$ in $G$, a contradiction. So a near leaf that
lies on such a path 
must be directly adjacent to $x_3$. Suppose first that
$d=4$. Since $q \geq 3$, and no near leaves are attached to
either $x_2$ or $x_d=x_4$, then $a=q-2 \geq 1$. Note that $G''$ is a
graph with $a+1$ connected components corresponding to the near leaves
adjacent to $x_3$ and to $(x_4x_5)$, which is the path of length
$d-3$. Notice also that $x_2$ is an isolated vertex of $G''$. Thus by
\cite[Lemma 6.2.7]{Rafael} 
$\depth(R/(I:x)) \geq a+1+1=q=\lceil{\frac{3q}{3}}\rceil \geq
\lceil{\frac{3+q}{3}}\rceil = \lceil{\frac{4+q-1}{3}}\rceil$ since $q
\geq 3$. 

Now suppose $d \geq 5$. Then the diameter of $G''$ is at least $d-3$
and $G''$ has at least $q-a-1$ near leaves in the connected component
containing $P_{d-2}$. Thus $\depth(R/(I:x_2)) \geq
\lceil{\frac{d-3+q-a-1-1}{3}}\rceil +a+1$ since $x_2$ is isolated and there
are $a$ additional connected components. If $a \geq 1$, then
$\lceil{\frac{d-3+q-a-1-1}{3}}\rceil +a+1 \geq
\lceil{\frac{d+q-1}{3}}\rceil$ as desired. Suppose $a=0$. If $G''$ has
$q$ near leaves, then $\depth(R/(I:x_2)) \geq
\lceil{\frac{d-3+q-1}{3}}\rceil+1=\lceil{\frac{d+q-1}{3}}\rceil$ as
desired. If $G''$ has $q-1$ near leaves, then $x_5$ cannot be an
additional near leaf. If $d=5$ then $x_5$ was already a near
leaf. Since $q \geq 3$ and $a=0$, there must be a near leaf on a path adjacent
to $x_4$ other than $P_{d-2}$. Since $x_4$ is not a leaf of $G''$, the
diameter of $G''$ is at least $d-2$. If $d\geq 6$ and $x_5$ is not a
near leaf of $G''$, then either $x_4$ is not a leaf, or there is a
non-leaf other than $x_6$ adjacent to $x_5$. In either case, the
diameter of $G''$ is at least $d-2$. Thus $\depth(R/(I:x_2)) \geq
\lceil{\frac{d-2+q-1-1}{3}}\rceil +1 =
\lceil{\frac{d+q-1}{3}}\rceil$. 

The result now follows from applying the Depth Lemma to the sequence
$$0 \rar R/(I:x_2) \rar R/I \rar R/(I,x_2)\rar 0.$$
\qed

\begin{Corollary}\label{bonus}
Suppose $G$ is a forest with $p$ connected components $G_1,\ldots
,G_p$, and $I=I(G)$. Let $d_i$ be the diameter of
$G_i$, let $d=\max {d_i}$, and let $q$ be the number of near leaves of
a component of diameter $d$. Then $\depth(R/I^t)\geq \max \{
\lceil{\frac{d-t+q}{3}}\rceil +p-1,p \}$ for all $t\geq 1$.
\end{Corollary}

\demo
As before, assume $d=d_1$ and $q$ is the number of near leaves of
$G_1$. When counting near leaves in this proof, only those in minors
of $G_1$ will be considered. To simplify notation, let $s=\max \{
\lceil{\frac{d-t+q}{3}}\rceil +p-1,p \}$. If $d \leq 3$, then $q \leq
2$ and so the result holds by 
Theorem ~\ref{forest} for small values of $d$ (and thus for small
values of $n$). When $t=1$, the result follows from Lemma
~\ref{treeBonus} and \cite[Lemma 6.2.7]{Rafael}. Fix the notation as in
Theorem ~\ref{forest} and 
assume $q \geq 3$ and $d\geq 4$. Then
$(I^t:x_1y) =I^{t-1}$ as before, and $\depth(R/(I^t:x_1y)) \geq \max \{
\lceil{\frac{d-t+1+q}{3}}\rceil +p-1,p \}\geq s$ by induction on $t$. Now
$(I^t,y)=(J^t,y)$ where $J$ is the ideal of $G'$. Note that the
diameter of $G'$ is at least $d-2$, $G'$ has at least $q-1$ near
leaves, and $x_1$ is an isolated vertex. Thus by induction on $n$,
$\depth(R/(I^t,y)) \geq \max \{ 
\lceil{\frac{d-2-t+q-1}{3}}\rceil +p-1,p \}+1 \geq s$ as before. 

As in Theorem ~\ref{forest} since $d \geq 4$, $r \geq 2$ and
$(I_r^t:y)=I_r^t$. Notice that the diameter of the graph corresponding
to $I_r$ is at least $d-3$. Note that because $x_r$ lies on a path of
maximal length and has distance two from a leaf on that path, any
path that connects $x_r$ to a leaf and does not 
contain any other vertex in $P_{d+1}$ must have length at most
two. Thus any near leaves on such a path must be directly adjacent to
$x_r$. Let $a$ be the number of near leaves adjacent to $x_r$ bur not
on $P_{d+1}$. Then the graph of $I_r$ has $p+a$ connected components, and the
number of near leaves of $I_r$ in the connected component containing
$P_{d-2}$ is at least $q-a-1$ when $d \geq 5$ since $y$ is also no
longer a near leaf. If $d=4$, the number of near leaves of $I_r$ is
$0$ and $a=q-2$ since every near leaf is adjacent to $x_r$, including
both that lie on $P_{d+1}$. So for $d=4$, the minor associated to
$I_r$ consists of $P_2$, together with $a$ additional connected
components and at least one isolated vertex. Thus $\depth(R_r/I_r^t) \geq 
1+ a + p-1+1 \geq q+p-1 \geq \lceil{\frac{q+3}{3}}\rceil +p-1
=\lceil{\frac{d+q-1}{3}}\rceil +p-1$ for $d=4$ and $q \geq 3$.

Assume $d=5$. Then either the connected component containing
$P_{d-2}$ contains two near leaves, and thus has diamenter at least
$3=d-2$, or $a \geq 1$. In the first case, by induction,
$$\depth(R_r/I_r^t) = \depth (R_r'/I_r^t) +1 \geq $$
$$\max \left\{ \left\lceil{\frac{d-2+q-1-t}{3}}\right\rceil +p-1,p\right\}+1 \geq s.$$ 
In the second case, $a \geq
1$, so 
$$\depth(R_r/I_r^t) = \depth (R_r'/I_r^t) +1 \geq $$
$$\max \left\{
\left\lceil{\frac{d-3+q-a-1-t}{3}}\right\rceil +p-1,p\right\}+a+1 \geq s.$$ 
If $d>5$, then $d-3 \geq 4$
and so either $P_{d-3}$ contains two near leaves, one of which was not
a near leaf of $I$, or $P_{d-3}$ is not maximal, and so the diameter
of $I_r$ is at least $d-2$. In the first case, the number of near
leaves is $q-a$ and the diameter is at least $d-3$, so
$$\depth(R_r/I_r^t) = \depth (R_r'/I_r^t) +1 \geq $$
$$\max \left\{ 
\left\lceil{\frac{d-3+q-a-t}{3}}\right\rceil +p-1,p\right\}+a+1 \geq s.$$ 
In the second case, the
number of near leaves is $q-a-1$ and the diameter is at least $d-2$,
so again 
$$\depth(R_r/I_r^t) = \depth (R_r'/I_r^t) +1 \geq $$
$$\max \left\{
\left\lceil{\frac{d-2+q-a-1-t}{3}}\right\rceil +p-1,p\right\}+a+1 \geq s.$$ 
Now by reverse induction on
$j$, using the proof from Theorem ~\ref{forest},
$\depth(R_1/(I_1^t:y))=\depth(R/((I^t:y),x_1)) 
\geq s$. 
As in the
theorem, the result now follows from Lemma ~\ref{sequencePair}.
\qed

\begin{Remark}\label{lowerBound}
{\rm Notice that Theorem ~\ref{forest} and Corollary ~\ref{bonus}
  provide a lower bound on where the stability of $\depth(R/I^t)$
  occurs for a tree. As noted in Lemma ~\ref{depth1}, the equality given in
  \cite[Proposition 3.3]{EH} or \cite[Theorem 1.2]{HH}, implies $\depth(R/I^t)
  = 1$ for all $t$ sufficiently large. In general, no bounds are known
  on how large $t$ must be to guarantee equality, although special
  cases are known. For 
  example, if $I$ is a complete graph, $\depth(R/I^t)=0$ for all
  $t\geq 2$ since $\m \in \Ass(R/I^t)$ by \cite{AJ}. This also follows
  from \cite[Corollary 3.4]{HH} noting that complete graphs are of
  the form $I_{n,2}$ in the notation used there. In the case of a
  tree, Theorem ~\ref{forest} shows that 
  $\depth(R/I^t) \geq 2$ for $t \leq d-2$, or in the case of Corollary
  ~\ref{bonus}, for $t \leq d+q-4$. Thus the depths of the powers do not
  stabilize until at least the $d-1^{st}$ power.
}
\end{Remark}

\begin{Example}
{\rm Consider the graph on 11 vertices with edges 
$$\{
  x_1x_2,x_2x_3,x_3x_4,x_4x_5, x_3x_6,x_6x_7,x_3x_8,x_8x_9,x_3x_{10},
  x_{10}x_{11} \}.$$ 
Then $d=4$ and $q=5$, so
  $\lceil{\frac{d+q-1}{3}}\rceil=3$, but it is easy to check
  using a computer program such as Macaulay 2 \cite{M2} that
  $\depth(R/I)=5$ in this example. Thus the bound given in Corollary
  ~\ref{bonus} is not necessarily sharp. It guarantees that $\depth(R/I^2) \geq
  \lceil{\frac{d+q-2}{3}}\rceil=3$ while the actual depth is again
  $5$. A careful reading of the proof shows that this is expected for
  $d=4$ and $q$ large. However, for $t$ large, the bound gains
  accuracy. For $t=5$, the bound and actual depth of $R/I^5$ are both
  $2$, and for $t=6$, both bound and actual depth are $1$. 

For larger $d$, the improved bound can be quite
accurate. For example, consider the graph on $9$ vertices with edges
$$\{ x_1x_2,x_2x_3,x_3x_4,x_4x_5,x_5x_6,x_6x_7,x_4x_8,x_8x_9\}.$$ 
Here $d=6$ and $q=3$. The
improved bound and the actual depth of $R/I^t$ agree for all powers $t
\not=3$, $t\leq 6$, as can be checked on Macaulay 2 \cite{M2}. In
particular, for $t=5$ the bound accurately predicts 
depth $2$, and for $t=6$ the bound, and the actual depth, become
one. Note that in both examples, this improved bound accurately
predicts where the depth will drop to one.}
\end{Example}

\section{Acknowledgements}

The author would like to thank J\"{u}rgen Herzog for asking the question that
sparked my initial interest in the depths of powers of edge
ideals. The author
would also like to thank Nate Dean for useful conversations regarding
grapth theory and Ray Heitmann and Rafael Villarreal for helpful
suggestions regarding an early version of the manuscript.

\end{document}